\documentclass[a4paper,11pt]{amsart}
\usepackage{graphicx}
 \usepackage{mathptmx}
\usepackage{amsmath}
\usepackage{amssymb}
\usepackage{enumitem}
\usepackage{xcolor}

\newmuskip\pFqmuskip

\newcommand*\pFq[6][8]{%
  \begingroup 
  \pFqmuskip=#1mu\relax
  \mathcode`=\string"8000
  \begingroup\lccode`\~=`\,
  \lowercase{\endgroup\let~}\pFqcomma
  F^{#2}_{#3}{\left(\genfrac..{0pt}{}{#4}{#5}\bigg|#6\right)}%
  \endgroup
}
\newcommand{\pFqcomma}{\mskip\pFqmuskip}

\newtheorem{theorem}{Theorem}

\begin{document}

\title{Higher-order degenerate harmonic numbers related to degenerate Riemann zeta function}

\author{Taekyun  Kim}
\address{Department of Mathematics, Kwangwoon University, Seoul 139-701, Republic of Korea}
\email{tkkim@kw.ac.kr}

\author{Dae San  Kim }
\address{Department of Mathematics, Sogang University, Seoul 121-742, Republic of Korea}
\email{dskim@sogang.ac.kr}


\begin{abstract}
Recently, Kim-Kim investigated the degenerate harmonic numbers and the degenerate hyperharmonic numbers as degenerate versions of the harmonic numbers and the hyperharmonic numbers, respectively. The aim of this paper is to study the higher-order degenerate harmonic numbers and the higher-order degenerate hyperharmonic numbers as higher-order versions for the degenerate harmonic numbers and the degenerate hyperharmonic numbers, respectively. In addition, we study the higher-order alternating degenerate hyperharmonic numbers as an `alternating version' of the higher-order degenerate hyperharmonic numbers. In more detail, we find generating functions of them, explicit expressions for them and some relations among them for those three kinds of numbers.
\end{abstract}

\maketitle

\markboth{\centerline{\scriptsize Higher-order degenerate harmonic numbers related to degenerate Riemann zeta function}}
{\centerline{\scriptsize Taekyun Kim and Dae San Kim}}

\section{Introduction}

Carlitz initiated (see [4]) the study of degenerate Bernoulli and degenerate Euler numbers and polynomials as degenerate versions of Bernoulli and Euler numbers and polynomials. Recently, this pioneering work regained interests of some mathematicians, various degenerate versions of quite a few special numbers and polynomials were investigated and some interesting arithmetical and combinatorial results were obtained along with some of their applications to other disciplines (see [9-13] and the references therein). For example, Kim-Kim investigated the degenerate harmonic numbers and the degenerate hyperharmonic numbers as degenerate versions of the harmonic numbers and the hyperharmonic numbers, respectively (see [10,11]).

The aim of this paper is to study the higher-order degenerate harmonic numbers and the higher-order degenerate hyperharmonic numbers which are higher-order versions of the degenerate harmonic numbers and the degenerate hyperharmonic numbers, respectively. In addition, we study the higher-order alternating degenerate hyperharmonic numbers as an alternating version of the higher-order degenerate hyperharmonic numbers. In more detail, we derive generating functions of them, explicit expressions for them and some relations among them for those three kinds of numbers.

The outline of this paper is as follows. In Section 1, we recall the degenerate exponentials and degenerate logarithms. We remind the reader of the polylogarithm and the degenerate polylogarithm. Then we recall the harmonic numbers and the hyperharmonic numbers. We recall the degenerate harmonic numbers and the degenerate hyperharmonic numbers. We also mention the binomial inversion. Section 2 is the main result of this paper. We define the degenerate Riemann zeta function as a degenerate version of the usual Riemann zeta function $\zeta_{\lambda}(s)$. In Theorem 2, we find an explicit expression of the degenerate harmonic number as a finite sum. In Theorem 3, $\zeta_{-\lambda}(2)$ is expressed as the integral over $(0,\infty)$ of the function $\frac{\log_{\lambda}(1+x)}{x(1+x)}$. Next, we define the degenerate harmonic numbers of order $\alpha$, denoted by $H_{n,\lambda}(\alpha)$, which reduce to the degenerate harmonic numbers $H_{n,\lambda}=H_{n,\lambda}(1)$, for $\alpha=1$. The generating function of those numbers are derived in Theorem 4. In Theorem 5, $\frac{1}{n}H_{n,-\lambda}$ is expressed as a finite sum involving $H_{m,\lambda}(2)$. We introduce the alternating degenerate harmonic numbers of order $\alpha$, denoted by $\overline{H}_{n,\lambda}(\alpha)$, and obtain the generating function of those numbers in Theorem 6. We derive a relation between the alternating degenerate harmonic numbers of order $\alpha-1$ and the degenerate harmonic numbers of order $\alpha$ in Theorem 7. We next define the degenerate hyperharmonic numbers of order $\alpha$, denoted by $H_{n,\lambda}^{(r)}(\alpha)$, which reduce to the degenerate hyperharmonic numbers $H_{n,\lambda}^{(r)}(1)=H_{n,\lambda}^{(r)}$, for $\alpha=1$. The generating function of those numbers is obtained in Theorem 8. Finally, we conclude our paper in Section 3. For the rest of this section, we recall the facts that are needed throughout this paper.\par

\vspace{0.1cm}

For any $\lambda\in\mathbb{R}$, the degenerate exponentials are defined by
\begin{equation}
e_{\lambda}^{x}(t)=\sum_{n=0}^{\infty}(x)_{n,\lambda}\frac{t^{n}}{n!},\quad (\mathrm{see}\ [12,13]),\label{1}
\end{equation}
where
\begin{equation}
(x)_{0,\lambda}=1,\ (x)_{n,\lambda}=x(x-\lambda)\cdots\big(x-(n-1)\lambda)\big),\ (n\ge 1).\label{2}
\end{equation}
As the inverse of $e_{\lambda}(t)=e_{\lambda}^{1}(t)$, the degenerate logarithm is defined by
\begin{equation}
	\log_{\lambda}(1+t)=\sum_{n=1}^{\infty}\frac{(1)_{n,1/\lambda}\lambda^{n-1}}{n!}t^{n},\quad (\mathrm{see}\ [9]). \label{3}
\end{equation}
Note that $\lim_{\lambda\rightarrow 0}\log_{\lambda}(1+t)=\log(1+t),\ \lim_{\lambda\rightarrow 0}e_{\lambda}(t)=e^{t}$. \par
It is well known that the harmonic numbers are defined by
\begin{equation}
	H_{0}=1,\quad H_{n}=1+\frac{1}{2}+\cdots+\frac{1}{n},\quad (n\ge 1),\quad (\mathrm{see}\ [1-5]). \label{4}
\end{equation}
From \eqref{4}, we note that
\begin{equation}
-\frac{\log(1-t)}{1-t}=\sum_{n=1}^{\infty}H_{n}t^{n},\quad (\mathrm{see}\ [1-16]). \label{5}
\end{equation} \par
For $s\in\mathbb{C}$, the polylogarithm is defined by
\begin{equation}
\mathrm{Li}_{s}(z)=\sum_{n=1}^{\infty}\frac{z^{n}}{n^{s}}=z+\frac{z^{2}}{2^{s}}+\frac{z^{3}}{3^{s}}+\cdots,\quad (\mathrm{see}\ [2-7,9-11). \label{6}
\end{equation}
Recently, Kim-Kim introduced the degenerate polylogarithm which is defined by
\begin{equation}
\mathrm{Li}_{k,\lambda}(x)=\sum_{n=1}^{\infty}\frac{(-\lambda)^{n-1}(1)_{n,1/\lambda}x^{n}}{(n-1)!n^{k}},\quad (|x|<1),\quad (\mathrm{see}\ [9]).\label{7}
\end{equation}
Note that $\mathrm{Li}_{1,\lambda}(x)=-\log_{\lambda}(1-x)$. \par
In 1996, Conway introduced the hyperharmonic numbers $H_{n}^{(r)}$ of order $r,\ (n,r\ge 0)$, which are given by
\begin{equation}
H_{0}^{(r)}=0,\quad (r\ge 0),\quad H_{n}^{(0)}=\frac{1}{n},\quad (n\ge 1),\quad H_{n}^{(r)}=\sum_{k=1}^{n}H_{k}^{(r-1)},\quad (n,r\ge 1). \label{8}
\end{equation}
From \eqref{8}, we note that
\begin{equation}
H_{n}^{(r)}=\binom{n+r-1}{n}\Big(H_{n+r-1}-H_{r-1}\Big),\quad (r\ge 1),\quad (\mathrm{see}\ [6]).\label{9}
\end{equation} \par
The degenerate harmonic numbers are defined by
\begin{equation}
H_{0,\lambda}=0,\quad H_{n,\lambda}=\sum_{k=1}^{n}\frac{(-\lambda)^{k-1}(1)_{k,1/\lambda}}{k!},\quad (n\ge 1),\quad (\mathrm{see}\ [10,12]). \label{10}
\end{equation}
Note that $\lim_{\lambda\rightarrow 0}H_{n,\lambda}=H_{n},\ (n\ge 1)$. \par
For $n\ge 0$, and $r \ge 1$, the degenerate hyperharmonic numbers are defined by (see\ [10])
\begin{equation}
H_{0,\lambda}^{(r)}=0,\quad (r\ge 1),\quad H_{n,\lambda}^{(1)}= H_{n,\lambda},\quad(n \ge 1),\quad H_{n,\lambda}^{(r)}=\sum_{k=1}^{n}H_{k,\lambda}^{(r-1)},\quad (n \ge 1, r\ge 2). \label{11}	
\end{equation}
Thus, by \eqref{11}, we get
\begin{displaymath}
	H_{n,\lambda}^{(r)}=\frac{(-1)^{r-1}}{\binom{\lambda-1}{r-1}}\binom{n+r-1}{n}\Big(H_{n+r-1,\lambda}-H_{r-1,\lambda}\Big),\quad (\mathrm{see}\ [10,11]).
\end{displaymath} \par
The next binomial inversion theorem is well known and will be used several times.
\begin{theorem}
For any integer $n \ge 0$, we have
\begin{equation*}
a_{n}=\sum_{k=0}^{n}(-1)^{n-k}\binom{n}{k}b_{k} \Longleftrightarrow b_{n}=\sum_{k=0}^{n}\binom{n}{k}a_{k}.
\end{equation*}
\end{theorem}

\section{Higher-order degenerate harmonic numbers}
First, we define the {\it{degenerate Riemann zeta function}} by
\begin{equation}
\zeta_{\lambda}(s)= \sum_{n=1}^{\infty}\frac{(-\lambda)^{n-1}(1)_{n,1/\lambda}}{n^{s}(n-1)!},\label{12}	
\end{equation}
where $s\in\mathbb{C}$ with $\mathrm{Re}(s)>1$.
Note that
\begin{displaymath}
	\lim_{\lambda\rightarrow 0}\zeta_{\lambda}(s)=\sum_{n=1}^{\infty}\frac{1}{n^{s}}=\zeta(s)
\end{displaymath}
is the Riemann zeta function. \par
From \eqref{10}, we note that
\begin{equation}
-\frac{1}{1-t}\log_{\lambda}(1-t)=\sum_{n=1}^{\infty}H_{n,\lambda}t^{n},\quad (|t|<1). \label{13}
\end{equation}
By \eqref{3} and \eqref{13}, we easily get
\begin{equation}
\frac{(1)_{n,1/\lambda}}{n!}\lambda^{n-1}=\sum_{k=1}^{n}(-1)^{n-k}\binom{n}{k}H_{k,-\lambda}.\label{14}
\end{equation}
The equation \eqref{14} is equivalent to
\begin{equation}
H_{n,-\lambda}=\sum_{k=1}^{n}\binom{n}{k}\frac{(1)_{k,1/\lambda}}{k!}\lambda^{k-1},\quad (n\ge 1).\label{15}	
\end{equation}
Note that
\begin{equation}
H_{n}=\lim_{\lambda\rightarrow 0}H_{n,-\lambda}=\sum_{k=1}^{n}\binom{n}{k}\frac{(-1)^{k-1}}{k}. \label{15-1}
\end{equation} \par
From \eqref{3}, we note that

\begin{align}
&\sum_{n=0}^{\infty}\frac{(-\lambda)^{n}(1)_{n+1,1/\lambda}}{(n+1)!}t^{n}=\frac{1}{t}\sum_{n=1}^{\infty}\frac{(1)_{n,1/\lambda}(-\lambda)^{n-1}}{n!}t^{n}=-\frac{1}{t}\log_{\lambda}(1-t) 	\label{16}\\
&=-\frac{1}{t}\bigg(\frac{-\log_{-\lambda}\big(1+\frac{t}{1-t}\big)}{\big(1+\frac{t}{1-t}\big)\big(1-t\big)}\bigg)=-\frac{1}{t(1-t)}\sum_{m=1}^{\infty}(-1)^{m}H_{m,-\lambda}\bigg(\frac{t}{1-t}\bigg)^{m} \nonumber \\
&=\sum_{m=0}^{\infty}(-1)^{m}H_{m+1,-\lambda}t^{m}(1-t)^{-m-2} \nonumber \\
&=\sum_{m=0}^{\infty}(-1)^{m}H_{m+1,-\lambda}t^{m}\sum_{k=0}^{\infty}\binom{m+k+1}{k+1}t^{k} \nonumber \\
&=\sum_{n=0}^{\infty}(n+1)\bigg(\sum_{m=0}^{n}(-1)^{m}\frac{H_{m+1,-\lambda}}{m+1}\binom{n}{m}\bigg)t^{n}. \nonumber
\end{align}
By comparing the coefficients on both sides of \eqref{16}, we get
\begin{equation}
\frac{(1)_{n+1,1/\lambda}\lambda^{n}}{(n+1)^{2}n!}=\sum_{m=0}^{n}(-1)^{n-m}\binom{n}{m}\frac{H_{m+1,-\lambda}}{m+1}. \label{18}	
\end{equation}
From \eqref{18} and by binomial inversion, we obtain the following theorem.
\begin{theorem}
For $n\ge 0$, we have
\begin{displaymath}
\frac{H_{n+1,-\lambda}}{n+1}=\sum_{m=0}^{n}\binom{n}{m}\frac{(1)_{m+1,1/\lambda}}{(m+1)^{2}m!}\lambda^{m},
\end{displaymath}
and
\begin{displaymath}
\frac{(1)_{n+1,1/\lambda}}{(n+1)^{2}n!}\lambda^{n}=\sum_{m=0}^{n}(-1)^{n-m}\binom{n}{m}\frac{H_{m+1,-\lambda}}{m+1}.
\end{displaymath}
\end{theorem}
From \eqref{13}, we have
\begin{align}
\frac{\log_{\lambda}(1+x)}{x(1+x)}&=-\frac{1}{x}\bigg(-\frac{\log_{\lambda}(1+x)}{1+x}\bigg)=\sum_{n=1}^{\infty}H_{n,\lambda}(-x)^{n-1}\label{20}\\
&=\sum_{n=0}^{\infty}H_{n+1,\lambda}(-x)^{n}. \nonumber	
\end{align}
By \eqref{18} and \eqref{20}, we get
\begin{align}
\int_{0}^{\infty}\frac{\log_{\lambda}(1+x)}{x(1+x)}dx&=\int_{0}^{1}\frac{\log_{\lambda}\big(1+\frac{x}{1-x}\big)}{\big(\frac{x}{1-x}\big)\big(1+\frac{x}{1-x}\big)}\frac{1}{(1-x)^{2}}dx \label{21} \\
&=\sum_{m=1}^{\infty}H_{m,\lambda}(-1)^{m-1}\int_{0}^{1}\bigg(\frac{x}{1-x}\bigg)^{m-1}\frac{1}{(1-x)^{2}}dx \nonumber \\
&=\sum_{m=0}^{\infty}H_{m+1,\lambda}(-1)^{m}\int_{0}^{1}x^{m}\frac{1}{(1-x)^{m+2}}dx \nonumber \\
&=\sum_{m=0}^{\infty}H_{m+1,\lambda}(-1)^{m}\sum_{k=0}^{\infty}\binom{k+m+1}{m+1}\int_{0}^{1}x^{m+k}dx \nonumber\\
&=\sum_{n=0}^{\infty}\sum_{m=0}^{n}H_{m+1,\lambda}(-1)^{m}\binom{n+1}{m+1}\int_{0}^{1}x^{n}dx \nonumber \\
&=\sum_{n=0}^{\infty}\sum_{m=0}^{n}\frac{H_{m+1,\lambda}}{m+1}\binom{n}{m}(-1)^{m} \nonumber \\
&=\sum_{n=0}^{\infty}\frac{(1)_{n+1,-1/\lambda}}{(n+1)^{2}n!}\lambda^{n}=\sum_{n=1}^{\infty}\frac{(1)_{n,-1/\lambda}\lambda^{n-1}}{n^{2}(n-1)!} \nonumber \\
&=\mathrm{Li}_{2,-\lambda}(1)=\zeta_{-\lambda}(2). \nonumber
\end{align}
Therefore, by \eqref{21}, we obtain the following theorem.
\begin{theorem}
The following identities hold true.
\begin{displaymath}
\int_{0}^{\infty}\frac{\log_{\lambda}(1+x)}{x(1+x)}dx=\mathrm{Li}_{2,-\lambda}(1)=\zeta_{-\lambda}(2).
\end{displaymath}
\end{theorem}
Now, we define the {\it{degenerate harmonic numbers of order $\alpha(\in\mathbb{N})$}} which are given by
\begin{equation}
H_{0,\lambda}(\alpha)=0,\quad H_{n,\lambda}(\alpha)=\sum_{k=1}^{n}\frac{(-\lambda)^{k-1}(1)_{k,1/\lambda}}{k^{\alpha}(k-1)!},\quad(n \ge 1).\label{22}
\end{equation}
Note that $\lim_{\lambda\rightarrow 0}H_{n,\lambda}(\alpha)=\sum_{k=1}^{n}\frac{1}{k^{\alpha}}=H_{n}(\alpha)$ and $H_{n,\lambda}(1)=H_{n,\lambda},\ (n\ge 1)$. Here we recall that $H_{n}(\alpha)$ are called the generalized harmonic numbers of order $\alpha$ (see [8]).\\
From \eqref{7} and \eqref{22}, we have
\begin{align}
\frac{1}{1-t}\mathrm{Li}_{\alpha,\lambda}(t)&=\sum_{l=0}^{\infty}t^{l}\sum_{k=1}^{\infty}\frac{(-\lambda)^{k-1}(1)_{k,1/\lambda}}{k^{\alpha}(k-1)!}t^{k}	\label{23} \\
&=\sum_{n=1}^{\infty}\bigg(\sum_{k=1}^{n}\frac{(-\lambda)^{k-1}(1)_{k,1/\lambda}}{k^{\alpha}(k-1)!}\bigg)t^{n} \nonumber \\
&=\sum_{n=1}^{\infty}H_{n,\lambda}(\alpha)t^{n}.\nonumber
\end{align}
Thus we obtain the following theorem from \eqref{23}.
\begin{theorem}
For $|t|<1$, we have
\begin{equation}
\frac{1}{1-t}\mathrm{Li}_{\alpha,\lambda}(t)=\sum_{n=1}^{\infty}H_{n,\lambda}(\alpha)t^{n}.\label{24}
\end{equation}
\end{theorem}
By using \eqref{7} and \eqref{15}, we have
\begin{align}
\mathrm{Li}_{2,\lambda}\bigg(\frac{t}{1+t}\bigg)&=\sum_{m=1}^{\infty}\frac{(-\lambda)^{m-1}(1)_{m,1/\lambda}}{m^{2}(m-1)!}\bigg(\frac{t}{1+t}\bigg)^{m} \label{25} \\
&=\sum_{m=1}^{\infty}\frac{(-\lambda)^{m-1}(1)_{m,1/\lambda}}{m(m-1)!}t^{m}\frac{1}{m}\sum_{k=0}^{\infty}\binom{m+k-1}{m-1}(-1)^{k}t^{k}\nonumber\\
&=\sum_{m=1}^{\infty}\frac{(-\lambda)^{m-1}(1)_{m,1/\lambda}}{m!}t^{m}\sum_{k=0}^{\infty}\binom{m+k}{m}\frac{(-1)^{k}}{m+k}t^{k} \nonumber \\
&=\sum_{n=1}^{\infty}\frac{1}{n}\bigg(\sum_{m=1}^{n}\frac{(-\lambda)^{m-1}(1)_{m,1/\lambda}}{m!}\binom{n}{m}(-1)^{n-m}\bigg)t^{n}\nonumber \\
&=\sum_{n=1}^{\infty}\frac{(-1)^{n-1}}{n}H_{n,-\lambda}t^{n}. \nonumber
\end{align}
From \eqref{25}, we have
\begin{align}
&\sum_{n=1}^{\infty}\frac{(-1)^{n-1}}{n}H_{n,-\lambda}t^{n}=\mathrm{Li}_{2,\lambda}\bigg(\frac{t}{1+t}\bigg)=\frac{\mathrm{Li}_{2,\lambda}\big(\frac{t}{1+t}\big)}{\big(1-\frac{t}{1+t}\big)(1+t)}\label{26} \\
&=\sum_{m=1}^{\infty}H_{m,\lambda}(2) t^{m}\bigg(\frac{1}{1+t}\bigg)^{m+1}=\sum_{m=1}^{\infty}H_{m,\lambda}(2)t^{m}\sum_{k=0}^{\infty}\binom{m+k}{m}(-1)^{k}t^{k}\nonumber \\
&=\sum_{n=1}^{\infty}\bigg(\sum_{m=1}^{n}(-1)^{n-m}\binom{n}{m}H_{m,\lambda}(2)\bigg)t^{n}. \nonumber
\end{align}
Thus, by comparing the coefficients on both sides of \eqref{26}, we get
\begin{equation}
(-1)^{n-1}\frac{1}{n}H_{n,-\lambda}=\sum_{m=1}^{n}(-1)^{n-m}\binom{n}{m}H_{m,\lambda}(2),\quad (n\in\mathbb{N}). \label{27}
\end{equation}
The equation \eqref{27} is equivalent to
\begin{equation}
H_{n,\lambda}(2)=\sum_{m=1}^{n}\binom{n}{m}(-1)^{m-1}\frac{H_{m,-\lambda}}{m}.\label{28}
\end{equation}
Therefore, by \eqref{27} and \eqref{28}, we obtain the following theorem.
\begin{theorem}
For $n\in\mathbb{N}$, we have
\begin{displaymath}
\frac{1}{n}H_{n,-\lambda}=\sum_{m=1}^{n}\binom{n}{m}(-1)^{m-1}H_{m,\lambda}(2),
\end{displaymath}
and
\begin{displaymath}
H_{n,\lambda}(2)=\sum_{m=1}^{n}\binom{n}{m}(-1)^{m-1}\frac{H_{m,-\lambda}}{m}.
\end{displaymath}
\end{theorem}
From \eqref{16} and \eqref{18}, we note that
\begin{align}
-\int_{0}^{1}\frac{1}{x}\log_{\lambda}(1-x)dx&=\sum_{n=0}^{\infty}\bigg(\sum_{m=0}^{n}(-1)^{m}\frac{H_{m+1,-\lambda}}{m+1}\binom{n}{m}\bigg)(n+1)\int_{0}^{1}x^{n}dx \label{29}	\\
&=\sum_{n=0}^{\infty}\frac{(1)_{n+1,1/\lambda}(-\lambda)^{n}}{(n+1)^{2}n!}=\sum_{n=1}^{\infty}\frac{(1)_{n,1/\lambda}(-\lambda)^{n-1}}{n^{2}(n-1)!}=\zeta_{\lambda}(2). \nonumber
\end{align}
Thus, by \eqref{29}, we get
\begin{equation}
-\int_{0}^{1}\frac{1}{x}\log_{\lambda}(1-x)dx=\zeta_{\lambda}(2). \label{30}
\end{equation} \par
Now, we consider the {\it{alternating degenerate harmonic numbers of order $\alpha(\in\mathbb{N})$}} given by
\begin{equation}
\overline{H}_{0,\lambda}(\alpha)=0,\quad \overline{H}_{n,\lambda}(\alpha)=\sum_{k=1}^{n}\binom{n}{k}\frac{\lambda^{k-1}(1)_{k,1/\lambda}}{k^{\alpha}(k-1)!},\quad (n \ge 1). \label{31}
\end{equation}
From \eqref{23} and \eqref{31}, we note that
\begin{align}
-\frac{\mathrm{Li}_{\alpha,\lambda}(-\frac{t}{1-t})}{1-t}
&=\sum_{k=1}^{\infty}\frac{\lambda^{k-1}(1)_{k,1/\lambda}}{(k-1)!k^{\alpha}}t^{k}\frac{1}{(1-t)^{k+1}}\label{32} \\
&=\sum_{k=1}^{\infty}\frac{\lambda^{k-1}(1)_{k,1/\lambda}}{(k-1)!k^{\alpha}}t^{k}\sum_{l=0}^{\infty}\binom{k+l}{k}t^{l} \nonumber \\
&=\sum_{n=1}^{\infty}\sum_{k=1}^{n}\binom{n}{k}\frac{\lambda^{k-1}(1)_{k,1/\lambda}}{(k-1)!k^{\alpha}} t^{n} \nonumber \\
&=\sum_{n=1}^{\infty}\overline{H}_{n,\lambda}(\alpha)t^{n}. \nonumber
\end{align}
Thus we get the next theorem from \eqref{32}
\begin{theorem}
For $|t|<1$, we have
\begin{equation}
-\frac{\mathrm{Li}_{\alpha,\lambda}(-\frac{t}{1-t})}{1-t}=\sum_{n=1}^{\infty}\overline{H}_{n,\lambda}(\alpha)t^{n}. \label{32-1}
\end{equation}
\end{theorem}
Note that $\lim_{\lambda\rightarrow 0}\overline{H}_{n,\lambda}(\alpha)=\sum_{k=1}^{n}\binom{n}{k}\frac{(-1)^{k-1}}{k^{\alpha}}=\overline{H}_{n}(\alpha)$, which may be called the {\it{alternating harmonic numbers of order $\alpha$}}. Here we observe from \eqref{15-1} that $\overline{H}_{n}(1)=H_{n}(1)=H_{n}$.\par
We observe that
\begin{align}
\mathrm{Li}_{\alpha,\lambda}\bigg(\frac{t}{1+t}\bigg)&=\sum_{m=1}^{\infty}\frac{(1)_{m,1/\lambda}(-\lambda)^{m-1}}{m^{\alpha}(m-1)!}t^{m}\bigg(\frac{1}{1+t}\bigg)^{m} \label{33} \\
& \sum_{m=1}^{\infty}\frac{(1)_{m,1/\lambda}(-\lambda)^{m-1}}{m^{\alpha-1}(m-1)!}t^{m}\sum_{k=0}^{\infty}\binom{m+k}{m}\frac{(-1)^{k}}{m+k}t^{k}\nonumber \\
&=\sum_{n=1}^{\infty}\frac{(-1)^{n-1}}{n}\bigg(\sum_{m=1}^{n}\frac{(1)_{m,1/\lambda}\lambda^{m-1}}{m^{\alpha-1}(m-1)!}\binom{n}{m}\bigg)t^{n} \nonumber \\
&=\sum_{n=1}^{\infty}(-1)^{n-1}\frac{\overline{H}_{n,\lambda}(\alpha-1)}{n}t^{n}. \nonumber
\end{align}
Thus, by \eqref{23} and \eqref{33}, we get
\begin{align}
\sum_{n=1}^{\infty}(-1)^{n-1}\frac{\overline{H}_{n,\lambda}(\alpha-1)}{n}t^{n}&=\frac{\mathrm{Li}_{\alpha,\lambda}\big(\frac{t}{1+t}\big)}{\big(1-\frac{t}{1+t}\big)(1+t)}=\sum_{m=1}^{\infty}H_{m,\lambda}(\alpha)t^{m}\bigg(\frac{1}{1+t}\bigg)^{m+1}\label{34} \\
&=\sum_{m=1}^{\infty}H_{m,\lambda}(\alpha)t^{m}\sum_{k=0}^{\infty}\binom{m+k}{m}(-1)^{k}t^{k} \nonumber \\
&=\sum_{n=1}^{\infty}\bigg(\sum_{m=1}^{n}H_{m,\lambda}(\alpha)(-1)^{n-m}\binom{n}{m}\bigg)t^{n}.\nonumber
\end{align}
By comparing the coefficients on both sides of \eqref{34}, we get
\begin{equation}
(-1)^{n-1}\frac{\overline{H}_{n,\lambda}(\alpha-1)}{n}=\sum_{m=1}^{n}H_{m,\lambda}(\alpha)(-1)^{n-m}\binom{n}{m},\quad (n\in\mathbb{N}). \label{35}	
\end{equation}
The equation \eqref{35} is equivalent to
\begin{equation}
H_{n,\lambda}(\alpha)=\sum_{m=1}^{n}\frac{\overline{H}_{m,\lambda}(\alpha-1)}{m}(-1)^{m-1}\binom{n}{m}. \label{36}
\end{equation}
Therefore, by \eqref{35} and \eqref{36}, we obtain the following theorem.
\begin{theorem}
For $\alpha\in\mathbb{N}$, we have
\begin{displaymath}
\frac{1}{n}\overline{H}_{n,\lambda}(\alpha-1)=\sum_{m=1}^{n}\binom{n}{m}(-1)^{m-1}H_{m,\lambda}(\alpha),
\end{displaymath}
and
\begin{displaymath}
H_{n,\lambda}(\alpha)=\sum_{m=1}^{n}\frac{\overline{H}_{m,\lambda}(\alpha-1)}{m}\binom{n}{m}(-1)^{m-1}.
\end{displaymath}
\end{theorem}
Now, we consider the {\it{degenerate hyperharmonic numbers of order $\alpha$}}, $H_{n,\lambda}^{(r)}(\alpha),\,\, (n\ge 0,\ \alpha,r\in\mathbb{N})$, which are defined by
\begin{equation}
H_{0,\lambda}^{(r)}(\alpha)=0,\,\,(r \ge 1),\,\, H_{n,\lambda}^{(1)}(\alpha)=H_{n,\lambda}(\alpha),\,\,(n \ge 1),\,\, H_{n,\lambda}^{(r)}(\alpha)=\sum_{k=1}^{n}H_{k,\lambda}^{(r-1)}(\alpha),\,\, (n \ge 1, r\ge 2). \label{37}	
\end{equation}
From \eqref{37}, we have
\begin{align}
\frac{\mathrm{Li}_{\alpha,\lambda}(t)}{(1-t)^{r}}&=\frac{1}{(1-t)^{r-1}}\frac{\mathrm{Li}_{\alpha,\lambda}(t)}{1-t}=\frac{1}{(1-t)^{r-1}}\sum_{n=1}^{\infty}H_{n,\lambda}(\alpha)t^{n} \label{38} \\
&=\frac{1}{(1-t)^{r-2}}\sum_{n=1}^{\infty}H_{n,\lambda}^{(2)}(\alpha)t^{n}	
=\cdots=\sum_{n=1}^{\infty}H_{n,\lambda}^{(r)}(\alpha)t^{n}.\nonumber
\end{align}
Therefore, by \eqref{38}, we obtain the following theorem.
\begin{theorem}
For $|t|<1$, we have
\begin{displaymath}
\frac{\mathrm{Li}_{\alpha,\lambda}(t)}{(1-t)^{r}}=\sum_{n=1}^{\infty}H_{n,\lambda}^{(r)}(\alpha)t^{n},\quad (r,\alpha\in\mathbb{N}).
\end{displaymath}
\end{theorem}
It is known that the degenerate hyperharmonic numbers have the following identity:
\begin{align}
H_{n,\lambda}^{(k+1)}&=\frac{(-1)^{k}}{\binom{\lambda-1}{k}}\binom{n+k}{k}\Big(H_{n+k,\lambda}-H_{k,\lambda}\Big)\label{39} \\
&=\frac{\binom{n+k}{k}}{\binom{k-\lambda}{k}}\Big(H_{n+k,\lambda}-H_{k,\lambda}\Big),\quad (\mathrm{see}\ [10,11]).\nonumber	
\end{align}
Here we may naturally ask whether it is possible to obtain such an identity for the generalized degenerate hyperharmonic numbers as the equation \eqref{39}.
So we leave the following as an open question to the researchers. \\

\vspace{0.1cm}

{ \bf Open Question}: Is it possible to get an expression of the  form
\begin{displaymath}
H_{n,\lambda}^{(r+1)}(\alpha)=A\Big(H_{n+k,\lambda}(\alpha)-H_{k,\lambda}(\alpha)\Big),\quad(\alpha,k\in\mathbb{N}),
\end{displaymath}
for some value of $A$?

\section{Conclusion}
In this paper, we defined the degenerate harmonic numbers of order $\alpha$, $H_{n,\lambda}(\alpha)$, as the sum of first $n$ terms of $\zeta_{\lambda}(\alpha)$, namely the degenerate Riemann zeta function evaluated at $\alpha$. Moreover, we introduced the alternating degenerate harmonic numbers of order $\alpha$ and the degenerate hyperharmonic numbers of order $\alpha$. For those three kinds of numbers, we investigated generating functions of them, explicit expressions for them and some relations among them. We suggested an open problem about expressing the degenerate hyperharmonic numbers of order $\alpha$ in terms of the degenerate harmonic numbers of order $\alpha$, which is possible for $\alpha=1$ (see \eqref{11}, \eqref{39}). \par
We would like to continue to explore various degenerate versions of many special polynomials and numbers with applications to physics, science and engineering as well as to mathematics.

\end{document}